\documentclass[letterpaper, 10 pt, journal, twoside]{ieeetran}

\usepackage{cite}
\usepackage{amsmath,amssymb,amsfonts}
\usepackage{algorithm,algorithmic}
\usepackage{mathabx}
\usepackage{bm}
\usepackage{graphicx}
\usepackage{textcomp}
\usepackage{color}

\newcommand{\hN}{\mathcal{N}}

\newcommand{\cL}{{\cal{L}}}
\newcommand{\cN}{{\cal N}}

\newcommand{\cE}{{\cal E}}

\newcommand{\cK}{{\cal K}}

\newtheorem{theorem}{Theorem}

\newtheorem{lemma}{Lemma}

\newtheorem{remark}{Remark}
\newtheorem{assumption}{Assumption}
\newtheorem{proposition}{Proposition}

\renewcommand{\baselinestretch}{0.98}

\begin{document}
\title{\LARGE Economic Dispatch With Distributed Energy Resources:\\
	\Large Co-Optimization of Transmission and Distribution Systems}
\author{Xinyang Zhou, \IEEEmembership{Member, IEEE}, Chin-Yao Chang, \IEEEmembership{Member, IEEE}, Andrey Bernstein, \IEEEmembership{Member, IEEE}, Changhong Zhao, \IEEEmembership{Member, IEEE}, and Lijun Chen, \IEEEmembership{Member, IEEE}
\thanks{This work was authored in part by the National Renewable Energy Laboratory, operated by Alliance for Sustainable Energy, LLC, for the U.S. Department of Energy (DOE) under Contract No. and DE-AC36-08GO28308 and DE-EE-0007998. Funding provided by U.S. Department of Energy Office of Electricity Delivery and Energy Reliability Advanced Grid Research \& Development through 2020 Maximizing Sensor Measurement Data through Adaptive Real Time Control project at NREL. Funding provided by U.S. Department of Energy Office of Energy Efficiency and Renewable Energy Solar Energy Technologies Office. The views expressed in the article do not necessarily represent the views of the DOE or the U.S. Government. The U.S. Government retains and the publisher, by accepting the article for publication, acknowledges that the U.S. Government retains a nonexclusive, paid-up, irrevocable, worldwide license to publish or reproduce the published form of this work, or allow others to do so, for U.S. Government purposes.}
\thanks{X. Zhou, C.-Y, Chang, and A. Bernstein are with the National Renewable Energy Laboratory, Golden, CO 80401, USA (Emails: \{xinyang.zhou, chinyao.chang, andrey.bernstein\}@nrel.gov).}
\thanks{C. Zhao is with the Department of Information Engineering, the Chinese University of Hong Kong, HKSAR, China (Email: chzhao@ie.cuhk.edu.hk).}
\thanks{L. Chen is with the College of Engineering and Applied Science, University of Colorado, Boulder, CO 80309, USA  (Email: lijun.chen@colorado.edu).}
}

\pagestyle{empty}

\maketitle
\thispagestyle{empty}

\begin{abstract}
The increasing penetration of distributed energy resources (DERs) in the distribution networks has turned the conventionally passive load buses into active buses that can provide grid services for the transmission system. To take advantage of the DERs in the distribution networks, this letter formulates a transmission-and-distribution (T\&D) systems co-optimization problem that achieves economic dispatch at the transmission level and optimal voltage regulation at the distribution level by leveraging large generators and DERs. A primal-dual gradient algorithm is proposed to solve this optimization problem jointly for T\&D systems, and a distributed market-based equivalent of the gradient algorithm is used for practical implementation. The results are corroborated by numerical examples with the IEEE 39-Bus system connected with 7 different distribution networks.
\end{abstract}


\section{Introduction}

The rising electricity demand and the shortage of power supply have caused surging electricity prices and even blackouts in peak hours; a few unfortunate events have occurred in recent years with or without market manipulation \cite{joskow2001california}. Meanwhile, the penetration of distributed energy resources (DERs) has been deepening in distribution systems, with residential photovoltaic (PV) devices, energy storage devices, and electric vehicles (EVs) becoming increasingly popular. Such DERs can potentially meet (part of) the demand from the distribution networks, and provide grid services such as voltage regulation. Involving residential DERs for energy supply without disturbing distribution system operation becomes operationally desired and economically sensible for the overall transmission-and-distribution (T\&D) systems.

In the literature, joint generator-side and load-side control has been proposed to assist power balancing and frequency regulation in the transmission systems \cite{mallada2017optimal, zhao2014design, wu2017hierarchical, kasis2016primary,chen2017reverse}. These works usually focus on dynamics in the transmission system by treating load buses as controllable nodes without detailing the distribution system structure at the load buses.

Optimizing DERs in distribution systems has been extensively studied in the past decade. Various problem formulations and solution methods have been proposed to optimally coordinate DERs for voltage regulation, loss minimization, dispatching signal tracking \cite{gan2014exact, vsulc2014optimal, molzahn2017survey}, etc. Works such as \cite{dall2018optimal,bernstein2019} propose a concept of virtual power plant that enables  the distribution network to provide a certain amount of aggregate power output by coordinating the DERs within the network. Most of these works usually focus on distribution system analysis and do not model any transmission structure.

There are a few works on T\&D co-optimization. In \cite{caramanis2016co},  concrete models for the transmission network, the distribution networks, and DERs are formulated, and  a multi-level solution method to solve the subproblems for each layer in sequence is proposed. In \cite{papavasiliou2018coordination}, a T\&D coordination scheme is proposed by solving respective subproblems for the two levels. However, there is no well-formulated joint T\&D optimization problem proposed in existing works, so it is difficult to characterize the global performance of their solutions. Moreover, solving subproblems for transmission and distribution networks in sequence might be suboptimal compared to the solution obtained from solving the joint T\&D co-optimization problem because the latter usually has a larger feasible set to find solutions.

In this letter, we first formulate a convex optimization problem featuring economic dispatch at the transmission level while ensuring optimal voltage regulation at the distribution level with linearized power flow equations. The outputs of the large generators in the transmission network and those of the DERs in the distribution networks are jointly optimized. Next, we propose a primal-dual gradient algorithm to solve the optimization problem with provable convergence. A market-based distributed implementation of the gradient algorithm is then designed to provide practical application. Finally, we illustrate the performance of the proposed scheme on the IEEE 39-bus system connected with 7 different distribution feeders.

The rest of this letter is structured as follows. Section~\ref{sec:model} models the T\&D networks. Section~\ref{sec:coopf} formulates the T\&D system co-optimization problem and proposes a gradient algorithm for solving it. Section~\ref{sec:market} designs a market-based distributed implementation of the gradient algorithm. Section~\ref{sec:num} presents numerical results and Section~\ref{sec:conclude} concludes this letter.

\emph{Notation:} We use bold uppercase letters to represent matrices, e.g., $\mathbf{A}$; italic bold letters to represent vectors, e.g, $\bm{A}$ and $\bm{a}$; and non-bold letters to represent scalars, e.g., $A$ and $a$. Superscript $^{\top}$ performs vector or matrix transpose. 
$[\cdot]_{\Omega}$ makes projection upon set $\Omega$. Operator $\bigtimes$ represents the Cartesian product of sets. 

\section{System Model}\label{sec:model}

In this section, we provide the model for both transmission and distribution systems, where the term ``bus" or ``control area" is used for the former, and ``node" is used for the latter. 

Consider a power transmission network, denoted by a graph $(\cK,\cE)$, where $\cK=\{1,\ldots,K\}$ is a set of buses or control areas; and the set $\cE\subset\cK\times\cK$ collects undirected transmission lines connecting the buses. 

Without loss of generality, for each bus $k\in\cK$, we assume there to be a dispatchable generator $k$ with mechanical power input $P_k^M$ and a distribution feeder indexed with $k$ with a total real power load $P^L_k$ injected at its substation. Define $\bm{P}^M:=[P_k^M]^{\top}_{k\in\cK}$ and $\bm{P}^L:=[P_k^L]^{\top}_{k\in\cK}$. Denote by $P_k^0$ the remaining uncontrollable power injection at bus $k$. We assume that the transmission system is lossless to have the following power balance equation:
\begin{eqnarray}
\sum_{k\in \cK} \Big(P_k^M - P_k^L(\bm{p}_k,\bm{q}_k)+P_k^0\Big) = 0.
\end{eqnarray}

Distribution feeder $k$ has a radial topology $(\cN_k,\cE_k)$ with a set $\cN_k$ collecting all its $N_k$ nodes  and a set $\cE_k$ collecting their connecting distribution lines. 
Let $\bm{v}_k:=[v_{k,1},\ldots,v_{k,N_k}]^{\top}\in\mathbb{R}^{N_k}$ denote the voltage magnitudes vector in the distribution system, and $\bm{p}_k:=[p_{k,1},\ldots,p_{k,N_k}]^{\top}\in\mathbb{R}^{N_k}$ and $\bm{q}_k:=[q_{k,1},\ldots,q_{k,N_k}]^{\top}\in\mathbb{R}^{N_k}$ the real and reactive power injections from all its nodes. For the purpose of algorithms design, we leverage a linear power flow model\footnote{The linear model is used only to formulate the optimization problem and devise an efficient solution algorithm. The  simulation experiments in Section~\ref{sec:num} are performed using the exact (AC) power flow model.}:
\begin{eqnarray}
\bm{v}_k &=& \mathbf{A}_k\bm{p}_k+\mathbf{B}_k\bm{q}_k+\bm{c}_k,\label{eq:plk}\\
P^L_k &=& \bm{M}_k^{\top}\bm{p}_k+\bm{N}_k^{\top}\bm{q}_k+d_k.\label{eq:x0}
\end{eqnarray}
where $\mathbf{A}_k,\mathbf{B}_k\in\mathbb{R}^{N_k\times N_k}$, $\bm{c}_k,\bm{M}_k,\bm{N}_k\in\mathbb{R}^{N_k}$, and $d_k\in\mathbb{R}$ are system parameters that can be computed using methods such as \cite{bolognani2016,bernstein2018}. In the following, we  use $\bm{v}_k(\bm{p}_k,\bm{q}_k)$ and $P_k^L(\bm{p}_k,\bm{q}_k)$ to represent \eqref{eq:plk} and \eqref{eq:x0}, respectively. We refer to Fig.~\ref{fig:IEEE39} for an illustrative system setup.

\section{Problem Formulation and Gradient Algorithm}\label{sec:coopf}

\subsection{Controllable Devices}
We assume that generator $k$ has a cost function $C_k^M(P_k^M)$ and a feasible set featuring its operational limits $\Omega^M_k$. Meanwhile, distribution feeder $k$ controls its total injected power $P_k^L$ indirectly through $\bm{p}_k$ and $\bm{q}_k$ from DERs while ensuring its voltage constraints modeled as $\bm{g}_k(\bm{v}_k(\bm{p}_k,\bm{q}_k)) \leq \bm{0}_{m_k}$. Similar to the generators, each DER $i$ has a cost function denoted by $C_{k,i}(p_{k,i},q_{k,i})$ and a feasible set $\Omega_{k,i}$. Define $\bm{\Omega}^M:=\bigtimes_{k\in\cK} \Omega^M_k$ and $\bm{\Omega}:=\bigtimes_{k\in\cK}\bigtimes_{i\in\cN_k} \Omega_{k,i}$. We have the following assumption on the cost and constraints functions.

\begin{assumption}	\label{ass:device}
	$\Omega^M_k$ is convex for all $k\in\cK$. $\Omega_{k,i}$ is convex for all $i\in\cN_k, k\in\cK$. Functions $C_{k,i}(p_{k,i},q_{k,i})$ for all $i\in\hN_k, k\in\cK$ and $C^M_k(P_k^M)$ for all $k\in\cK$ are continuously differentiable and strongly convex in $(p_{k,i},q_{k,i})$ and $P_k^M$, respectively, with bounded first-order derivatives. Functions $\bm{g}_k(\bm{v}_k)$ are convex and differentiable functions of $\bm{v}_k$ for all $k\in\cK$.
\end{assumption}

\subsection{T\&D Co-Optimization Problem}
In this part, we formulate a T\&D co-optimization problem that achieves economic dispatch over all generators and DERs while maintaining voltage constraints by DERs in the distribution feeders.

Let $N=\sum_{k\in\cK}N_k$ be the total number of nodes in all distribution networks. Denote by  $\bm{p}=[\bm{p}_1^{\top},\ldots,\bm{p}_K^{\top}]^{\top},\ \bm{q}=[\bm{q}_1^{\top},\ldots,\bm{q}_K^{\top}]^{\top}\in\mathbb{R}^N$. Consider the following optimization problem subject to power flow and operational constraints:
\begin{subequations}\label{eq:R1}
	\begin{eqnarray}
	\hspace{-10mm}&\min& \sum_{k\in\cK}\hspace{-.5mm}\Big(\hspace{0mm}\sum_{i\in\cN_k} \hspace{-1.5mm}C_{k,i}(p_{k,i},q_{k,i})+C^M_k(P_k^M)\Big),\label{eq:obj1}\\
	\hspace{-10mm}&\text{over}&(\bm{p},\bm{q})\in\bm{\Omega},\bm{P}^M\in\bm{\Omega}^M\nonumber\\
	\hspace{-10mm}&\text{s.t.}&  \sum_{k\in \cK} \Big(P_k^M - P_k^L(\bm{p}_k,\bm{q}_k)+P_k^0\Big) = 0,\label{eq:balance1}\\[-2pt]
	\hspace{-10mm}&&\bm{g}_k(\bm{v}_k(\bm{p}_k,\bm{q}_k)) \leq \bm{0}_{m_k},\ \forall k\in\cK, \label{eq:volt1}
	\end{eqnarray}
\end{subequations}
where the cost function~\eqref{eq:obj1} adds up the generation costs of all large generators and small DERs in the distribution feeders, the equality constraint~\eqref{eq:balance1} ensures that power demand and supply are balanced, and the inequality constraint \eqref{eq:volt1} confines voltage magnitudes to within acceptable ranges.

Let $m=\sum_{k\in\cK}m_k$ be the dimension of distribution network constraints. Introduce dual variables $\lambda\in\mathbb{R}$ for the equality constraint~\eqref{eq:balance1} and nonnegative vector $\bm{\mu}=[\bm{\mu}_k^{\top}]^{\top}_{k\in\cK}\in\mathbb{R}_+^{m}$ for the inequality constraints~\eqref{eq:volt1} to have the following \emph{regularized} Lagrangian of \eqref{eq:R1} with a small constant $\eta>0$:
\begin{eqnarray}
\hspace{-3mm}&&\hspace{-3mm}\cL(\bm{p},\bm{q},\bm{P}^M;\lambda, \bm{\mu})\nonumber\\
\hspace{-3mm}&=&\hspace{-3mm}\sum_{k\in\cK}\bigg(\sum_{i\in\cN_k} C_{k,i}(p_{k,i},q_{k,i})+C^M_k\big(P_k^M\big)+ \nonumber\\[-3pt]
\hspace{-3mm}&&\hspace{-3mm} \bm{\mu}_k^{\top}\bm{g}_k\big(\bm{v}_k(\bm{p}_k,\bm{q}_k)\big)\Big) + \lambda \big(\sum_{k\in \cK} P_k^M - P_k^L(\bm{p}_k,\bm{q}_k)+P_k^0\big)\nonumber\\[-5pt]
\hspace{-3mm}&&\hspace{-3mm} - \frac{\eta(\lambda^2+\|\bm{\mu}\|^2_2)}{2}.  \label{eq:L}
\end{eqnarray}

Introducing the regularization term $-\eta(\lambda^2+\|\bm{\mu}\|^2_2)/2$ ensures strong concavity of $\cL(\bm{p},\bm{q},\bm{P}^M;\lambda, \bm{\mu})$ with respect to the dual variables, as well as provable convergence of gradient-based algorithms with a constant stepsize. However, a discrepancy proportional to $\eta$ is also brought in, which can be negligible if $\eta$ is small. We refer to Proposition 3.1 of \cite{koshal2011multiuser} for detailed analytical characterization of the discrepancy. Note that $\eta =0$ is used in Section~\ref{sec:num}, and the numerical results converge well.

\subsection{Gradient-Based Algorithm Design}\label{sec:alg}
We next design a primal-dual gradient algorithm to solve for the unique saddle point of  \eqref{eq:L}. For notational simplicity, we let $\bm{x}=[\bm{p}^{\top},\bm{q}^{\top},(\bm{P}^{M})^{\top}]^{\top}$ collect all the primal variables and $\bm{y}=[\lambda, \bm{\mu}^{\top}]^{\top}$ collect all the dual variables. Then, the iterative primal-dual gradient algorithm for finding the unique saddle point of the regularized Lagrangian~\eqref{eq:L} is given by:
\begin{subequations}\label{eq:primaldual}
	\begin{eqnarray}
	\bm{x}(t+1) & = & \left[\bm{x}(t) - \epsilon \frac{\partial\cL(\bm{x}(t);\bm{y}(t))}{\partial \bm{x}(t)}\right]_{\bm{\Omega}\times\bm{\Omega}^M},\\
	\bm{y}(t+1) & = & \left[\bm{y}(t) + \epsilon \frac{\partial\cL(\bm{x}(t);\bm{y}(t))}{\partial \bm{y}(t)}\right]_{\mathbb{R}\times\mathbb{R}_+^{m}}\hspace{-1mm},
	\end{eqnarray}
\end{subequations}
where $\epsilon$ is a constant stepsize and $t$ is the iteration index. The partial gradient of $\cL$ with respect to the decision variables are calculated as follows:
\begin{subequations}\label{eq:gradient}
	\begin{eqnarray}
	&\frac{\partial \cL}{\partial\bm{p}_k}=& \nabla_{\bm{p}_k} \hspace{-2mm} \sum_{i\in\cN_k} C_{k,i}(p_{k,i},q_{k,i})\nonumber\\[-2pt]
	&& - \lambda \bm{M}_k+\mathbf{A}^{\top}_k\nabla_{\bm{v}_k} \bm{g}_k(\bm{v}_k)^{\top}\bm{\mu}_k,\label{eq:gradientP}\\[2pt]
	&\frac{\partial \cL}{\partial\bm{q}_k}=&\nabla_{\bm{q}_k}\hspace{-2mm}\sum_{i\in\cN_k} C_{k,i}(p_{k,i},q_{k,i})\nonumber\\[-2pt]
	&& - \lambda \bm{N}_k + \mathbf{B}_k^{\top}\nabla_{\bm{v}_k} \bm{g}_k(\bm{v}_k)^{\top}\bm{\mu}_k,\label{eq:gradientQ}\\[2pt]
	&\frac{\partial \cL}{\partial P_k^M}=&d C_k^M(P_k^{M})/d {P_k^M}+\lambda, \label{eq:gradientPM}\\
	&\frac{\partial \cL}{\partial\bm{\mu}_k} =& \bm{g}_k\big(\bm{v}_k(\bm{p}_k,\bm{q}_k)\big) -\eta\bm{\mu}_k,\label{eq:gradientMu}\\
	&\frac{\partial \cL}{\partial \lambda} =&\sum_{k\in \cK} \big(P_k^M - P_k^L(\bm{p}_k,\bm{q}_k)+P_k^0\big) - \eta\lambda,\label{eq:gradientLambda}
	\end{eqnarray}
\end{subequations}
where $\nabla_{\bm{v}_k} \bm{g}_k(\bm{v}_k)$ is the Jacobian matrix of $\bm{g}_k$ with respect to $\bm{v}_k$, and Eqs.~\eqref{eq:gradientP}--\eqref{eq:gradientMu} are for all $k\in\cK$.

\subsection{Convergence Analysis}
Define a gradient operator $T(\bm{x};\bm{y})=\begin{bmatrix}{\partial\cL(\bm{x};\bm{y})}/{\partial \bm{x}} \\ -{\partial\cL(\bm{x};\bm{y})}/{\partial \bm{y}}\end{bmatrix}$. Based on Assumption~\ref{ass:device} and the regularization terms we have added for the dual variables to the Lagrangian, the next lemma follows by definition.
\begin{lemma}
	Based on Assumption~\ref{ass:device}, $T(\bm{x};\bm{y})$ is an $s$-strongly monotone operator with some constant $s>0$ and is $l$-Lipschitz continuous with some constant $l>0$ for any feasible $\bm{x}\in\bm{\Omega}\times\bm{\Omega}^M$ and $\bm{y}\in\mathbb{R}\times\mathbb{R}_+^{m}$.
\end{lemma}

\begin{theorem}\label{the:converge}
	Based on Assumption~\ref{ass:device}, given a constant stepsize $\epsilon$ such that $0 < \epsilon \leq \bar{\epsilon} < 2s/l^2$, the primal-dual gradient dynamics~\eqref{eq:primaldual} asymptotically converge to the unique saddle point of the regularized Lagrangian~\eqref{eq:L}.
\end{theorem}
We omit the detailed proof of Theorem~\ref{the:converge} here because of the space limit. It can be found in numerous literature, e.g., \cite{bertsekas1989parallel,boyd2004convex}. Moreover, asynchronous implementation of the proposed algorithm caused by communication delay or loss can be shown to converge to the same solutions under reasonable assumptions; we refer to \cite{bertsekas1989parallel} for more details.

\section{Distributed Market-Based Implementation}\label{sec:market}

\subsection{Economic Model}
Unlike most utility-owned power plants, the user-owned DERs in the distribution feeders are usually not obliged to follow the gradient steps specified in \eqref{eq:gradientP}--\eqref{eq:gradientQ}. On the contrary, users are naturally driven to minimize their own overall cost (or maximize their overall utility) featuring a trade-off between their DERs generation cost and electricity bills, which is formulated as follows:
\begin{eqnarray}
\hspace{-7mm} & \underset{\text{over}\  (p_{k,i},q_{k,i})\in\Omega_{k,i}}{\min f_{k,i}(p_{k,i},q_{k,i})}& \hspace{-3mm}  =C_{k,i}(p_{k,i},q_{k,i})+\alpha_{k,i}p_{k,i}+\beta_{k,i}q_{k,i}.\label{eq:userki}
\end{eqnarray}

Here, $\alpha_{k,i},\beta_{k,i}\in\mathbb{R}$ are the incentive signals/electricity prices for real and reactive power, respectively, set by the network operator for user $i$ in distribution feeder $k$. When $\alpha_{k,i}$ and $\beta_{k,i}$ are positive (resp. negative), users are incentivized to reduce (resp. increase) the values of their $p_{k,i}$ and $q_{k,i}$. Moreover, once the problem formulation---specifically, $C_{k,i}$ and $\Omega_{k,i}$---is revealed, $\alpha_{k,i}$ and $\beta_{k,i}$ can be used to induce certain values of $p_{k,i}$ and  $q_{k,i}$ by solving \eqref{eq:userki}.

However, private information of users is usually inaccessible to the network operator. We next present an iterative method to find the optimal signals to incentivize the users to react in a certain way that concurrently solves the T\&D optimization problem \eqref{eq:R1} \cite{zhou2017incentive}. As we will see, the resultant design can be seen as a market-based equivalent implementation of the primal-dual gradient algorithm~\eqref{eq:primaldual}.

\begin{figure}[t]
	\centering
	\includegraphics[width=.38\textwidth]{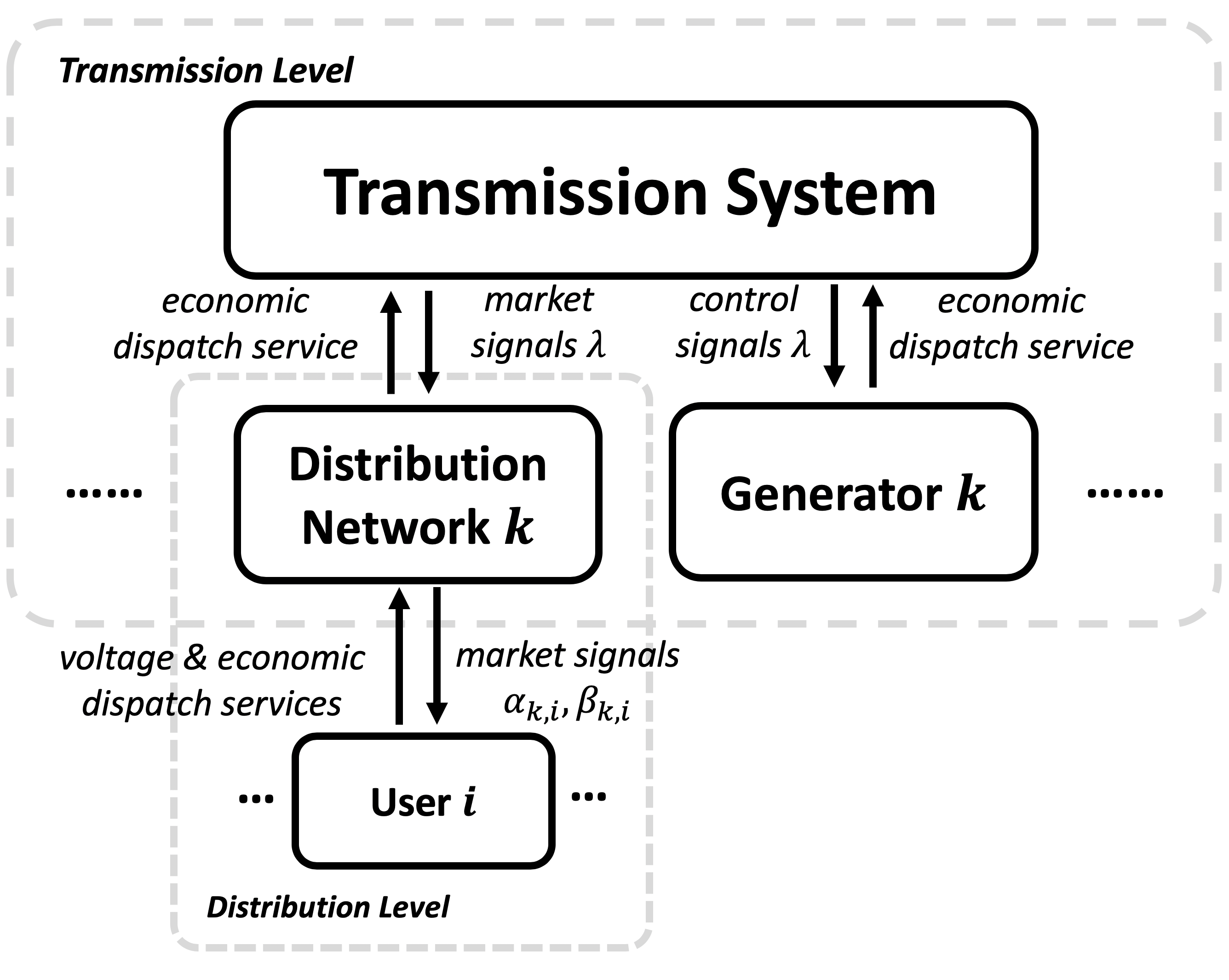}
	\caption{ Market-based distributed implementation of the primal-dual gradient algorithm for solving the T\&D system co-optimization problem.}
	\label{fig:flowchart}
	\vspace{-5mm}
\end{figure}

\subsection{Market-Based Distributed Implementation}
To incentivize the users to act according to \eqref{eq:primaldual} so that problem~\eqref{eq:R1} can be solved, network operator needs to carefully design the incentive signals. Note that when user $k,i$ solves its cost minimization problem~\eqref{eq:userki}, the gradient he/she takes is in the form of:
\begin{subequations}\label{eq:indigrad}
	\begin{eqnarray}
	\partial f_{k,i}/\partial p_{k,i}&=&\partial C_{k,i}/\partial p_{k,i} + \alpha_{k,i},\\
	\partial f_{k,i}/\partial q_{k,i}&=&\partial C_{k,i}/\partial q_{k,i} + \beta_{k,i}.
	\end{eqnarray}
\end{subequations}

Denote the signals vector of distribution feeder $k$ as $\bm{\alpha}_k = [\alpha_{k,1},\ldots,\alpha_{k,N_k}]^{\top}$ and $\bm{\beta}_k = [\beta_{k,1},\ldots,\beta_{k,N_k}]^{\top}$. By comparing \eqref{eq:indigrad} with \eqref{eq:gradientP}--\eqref{eq:gradientQ}, we design the incentive signals as:
\begin{subequations}\label{eq:signals}
	\begin{eqnarray}
	\bm{\alpha}_k&=& -\lambda \bm{M}_k +\mathbf{A}_k^{\top}\nabla_{\bm{v}_k} \bm{g}_k(\bm{v}_k)^{\top}\bm{\mu}_k,\label{eq:alpha}\\
	\bm{\beta}_k&=& -\lambda \bm{N}_k +\mathbf{B}_k^{\top}\nabla_{\bm{v}_k} \bm{g}_k(\bm{v}_k)^{\top}\bm{\mu}_k,\label{eq:beta}
	\end{eqnarray}
\end{subequations}
which relies on network information without any private information from the users. Using \eqref{eq:signals}, we propose a market-based iterative distributed algorithm presented as Algorithm~\ref{alg:distalg}. We illustrate the algorithm in Fig.~\ref{fig:flowchart} which also indicates the possibility of parallel execution of Algorithm~\ref{alg:distalg} for independent steps. By design, we have the following formal statement.
\begin{proposition}\label{prop1}
	Algorithm~\ref{alg:distalg} is equivalent to the primal-dual gradient algorithm for solving the saddle point of~\eqref{eq:L}.
\end{proposition}
Because Algorithm~\ref{alg:distalg} and the primal-dual gradient dynamics~\eqref{eq:primaldual} are equivalent, they share the same convergence properties. 

\begin{remark}
	For ease of presentation, we assume that all DERs update with gradient steps. In reality, however, non-cooperative DERs may update their setpoints by directly solving~\eqref{eq:userki} given the current incentive signals. This constitutes a dual ascend algorithm that is consistent with the market-based design~\eqref{eq:signals}. We refer to \cite{zhou2019online} for more details.
\end{remark}

\begin{algorithm}
	\caption{Distributed Market-Based T\&D Co-Optimization} \label{alg:distalg}
	\begin{algorithmic}
		
		\WHILE {stopping criterion not met}
		
		\STATE[S1] Given incentive signals $\alpha_{k,i}(t)$ and $\beta_{k,i}(t)$, user $i\in\cN_k, k\in\cK$ takes a gradient step toward solving his/her own optimization problem~\eqref{eq:userki} by:
		\begin{subequations}
			\begin{eqnarray*}
				p_{k,i}(t+1)\hspace{-3mm}&=&\hspace{-3mm} \Big[p_{k,i}(t)-\epsilon\big({\partial C_{k,i}(p_{k,i}(t),q_{k,i}(t))}/{\partial_{p_{k,i}}}\\[-4pt]
				&& + \alpha_{k,i}(t)\big)\Big]_{\Omega_{k,i}}, \\
				q_{k,i}(t+1)\hspace{-3mm}&=&\hspace{-3mm} \Big[q_{k,i}(t)-\epsilon\big({\partial C_{k,i}(p_{k,i}(t),q_{k,i}(t))}/{\partial_{q_{k,i}}}\\[-4pt]
				&& + \beta_{k,i}(t)\big)\Big]_{\Omega_{k,i}}. 
			\end{eqnarray*}
		\end{subequations}
		\STATE[S2] Dispatchable generator $k\in\cK$ updates its power setpoints by:
		\begin{eqnarray*}
			P_k^M(t+1)\hspace{-2mm}&=&\hspace{-2mm}\Big[P_k^M(t)-\epsilon\big(\frac{d C_k^M(P_k^{M}(t))}{d P_k^M}+\lambda(t)\big)\Big]_{\Omega_k^M}.
		\end{eqnarray*}
		
		\STATE [S3] Network operator updates the power flow by:
		\begin{subequations}
			\begin{eqnarray*}
				\hspace{-3mm}P^L_k(t+1) \hspace{-3mm}&=&\hspace{-3mm} \bm{M}_k^{\top}\bm{p}_k(t+1)+\bm{N}_k^{\top}\bm{q}_k(t+1)+d_k,\\
				\hspace{-3mm}\bm{v}_k(t+1) \hspace{-3mm}&=&\hspace{-3mm} \mathbf{A}_k\bm{p}_k(t+1)+\mathbf{B}_k\bm{q}_k(t+1)+\bm{c}_k.
			\end{eqnarray*}
		\end{subequations}
		
		\STATE[S4] Network operator updates the dual variables as:
		\begin{eqnarray*}
			\hspace{-5mm}\lambda(t+1)\hspace{-3mm}&=&\hspace{-3mm}\lambda(t)+\epsilon\Big(\sum_{k\in \cK} \big(P_k^M(t+1) -s P_k^L(t+1)+P_k^0\big) \\[-3pt]
			&&- \eta\lambda(t)\Big),\\
			\hspace{-5mm}\bm{\mu}_k(t+1) \hspace{-3mm}&=& \hspace{-3mm} \big[\bm{\mu}_k(t)+\epsilon\big(\bm{g}_k\big(\bm{v}_k(t+1)\big) -\eta\bm{\mu}_k(t)\big)\big]_{\mathbb{R}_+^{N_k}},
		\end{eqnarray*}
		and the incentive signals as:
		\begin{subequations}
			\begin{eqnarray*}
				\hspace{-5mm}\bm{\alpha}_k(t+1)\hspace{-3mm} &=& \hspace{-3mm}-\lambda(t+1) \bm{M}_k +\mathbf{A}_k^{\top}\nabla \bm{g}_k(\bm{v}_k(t+1))^{\top}\bm{\mu}_k(t+1),\label{eq:alphaupdate}\\
				\hspace{-5mm}\bm{\beta}_k(t+1) \hspace{-3mm}&=&\hspace{-3mm} -\lambda(t+1) \bm{N}_k +\mathbf{B}_k^{\top}\nabla \bm{g}_k(\bm{v}_k(t+1))^{\top}\bm{\mu}_k(t+1),\label{eq:betaupdate}
			\end{eqnarray*}
		\end{subequations}
		and sends the updated signals to the users.
		\ENDWHILE
	\end{algorithmic}
\end{algorithm}

\subsection{Nonlinear (AC) Power Flow Feedback}
In the primal-dual gradient update~\eqref{eq:primaldual}, as well as in Algorithm~\ref{alg:distalg}, while the Jacobian matrices ${\partial\cL(\bm{x};\bm{y})}/{\partial \bm{y}}$ and ${\partial\cL(\bm{x};\bm{y})}/{\partial \bm{x}}$ are calculated based on the linearized relationships~\eqref{eq:plk}--\eqref{eq:x0} to simplify the computation, they lose accuracy. To improve this, a feedback mechanism can be applied to reduce such modeling discrepancies and guarantee that the key system states such as $\bm{v}_k$ and $P_k^L$ are accurate. Specifically, when performing Algorithm~\ref{alg:distalg}, instead of using linearized power flow model in [S3], nonlinear power flow is used to calculate $\bm{v}_k$ and $P_k^L$, which are further fed into [S4] to update the dual variables. Stability analysis with the feedback mechanism can be found in \cite{zhou2017incentive}. We also use nonlinear (AC) power flow for the numerical results in Section~\ref{sec:num} next.

\section{Numerical Results}\label{sec:num}

\begin{figure}[t]
	\centering
	\includegraphics[width=.48\textwidth]{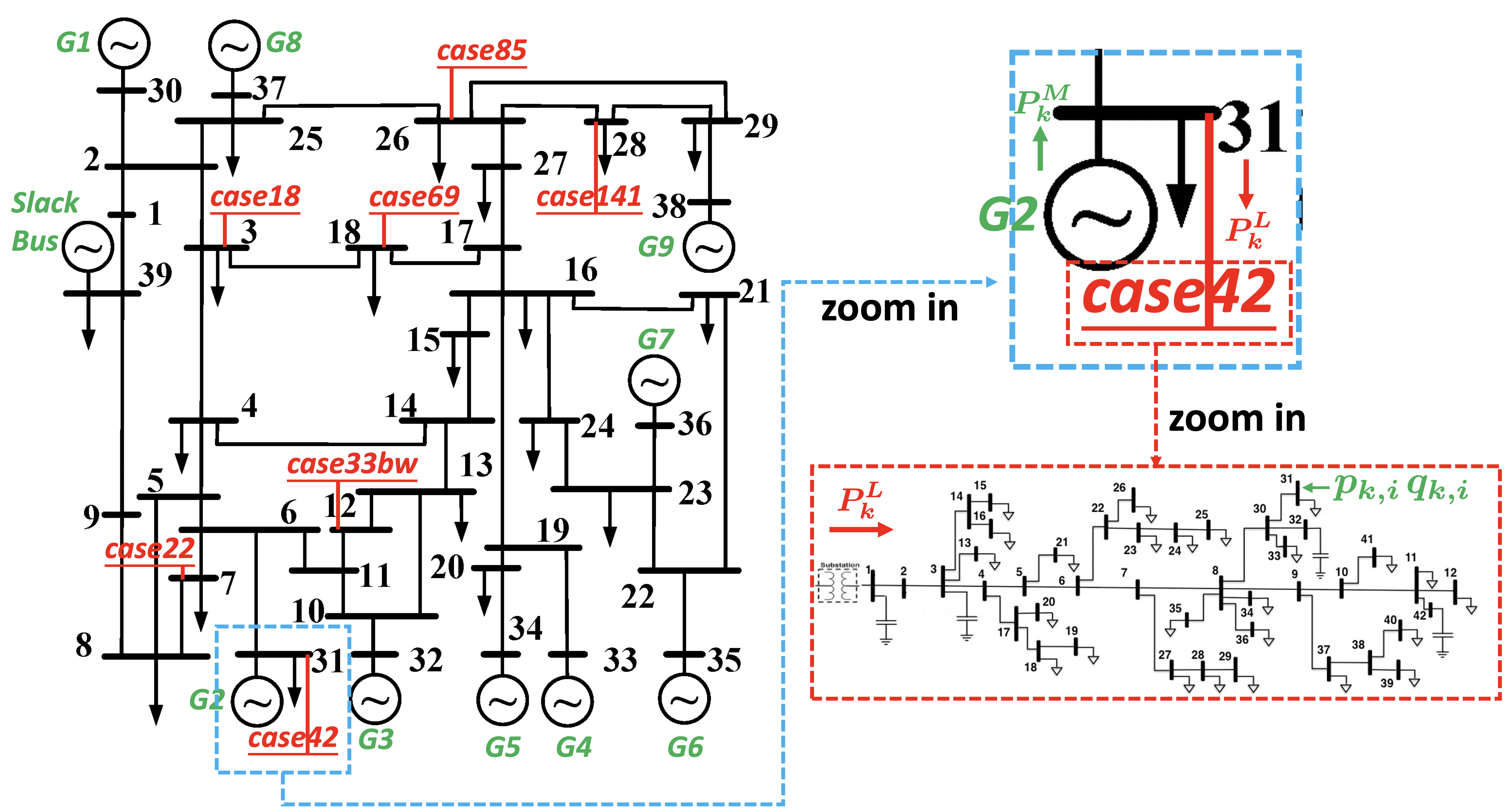}
	\caption{IEEE 39-Bus transmission system is used with 7 distribution networks connected to its load buses marked by red letters and 9 controllable generators connected to Buses 30--38. Bus 39 is the slack bus. Zoomed-in figures show the coupling between the transmission and distribution systems.}
	\label{fig:IEEE39}
	\vspace{-5mm}
\end{figure}

\subsection{System Setup}

\subsubsection{T\&D Systems}
We use the New England IEEE 39-Bus system as the transmission network with 9 controllable generators located at buses 30--38, and bus 39 (the slack bus) connected to the rest of US/Canada grid. We connect 7 different distribution networks---case18, case22, case33bw, case69, case85, case141 \cite{zimmerman2016matpower}, and SCE 42-bus system \cite{zhou2020reverse}---to load buses indexed 3, 7, 12, 18, 26, 28, and 31, respectively. The total power injected into the distribution feeders will be used as the demand of the corresponding load bus of the transmission system. See Fig.~\ref{fig:IEEE39} for the locations of the generators and distribution networks.

\subsubsection{Parameters}
We assign quadratic cost functions $c_k {P_k^{M}}^2$, with $c_k$ set to $1, 1.5, 1.3, 1.7, 1.8, 1, 2, 0.8$, and $1.2$, respectively, for the generators 1--9. We assign homogeneous cost functions $ p_{k,i}^2 + 0.1 q_{k,i}^2$ for DERs in distribution feeders.\footnote{Here, the reactive power cost can be interpreted as an opportunity cost because part of the DER's capacity is occupied.} The inequality constraints $\bm{g}_k(\bm{v}_k)\leq 0$ are set to $0.95~\text{p.u.}\leq \bm{v}_k \leq  1.05~\text{p.u.}$ for all $k$. Linearization parameters are based on LinDistFlow model to generate parameters $\mathbf{A}_k, \mathbf{B}_k, \bm{c}_k$ and $\bm{M}_k, \bm{N}_k, d_k$ in Eqs.~\eqref{eq:plk}--\eqref{eq:x0} for the distribution networks. Nonlinear power flow is solved to update both the transmission and the distribution systems every iteration with MATPOWER 7.0 \cite{zimmerman2016matpower}. We set $\eta=0$ and $\epsilon\approx 5\times 10^{-4}$, which is further slightly tuned for different feeders to improve convergence.

\subsubsection{Simulation Process} The system is initialized at a non-optimal point with voltage constraints violated in some distribution network. The system then approaches the optimal with primal-dual gradient dynamics before generator 7 is set to go down at iteration 10,000. At this point, we enlarge the feasible set of DERs by two times to allow them to contribute more to the grids without loss of generality.\footnote{Otherwise, DERs have already reached their previous limits because they are set to contribute cheaper power than the large generators.} The system will then approach a new optimal point with all constraints satisfied.

\subsubsection{Handling the Slack Bus}
In MATPOWER and other power system analysis tools, slack buses are essential to compensate loss and balance power; however, this automatic balancing functionality causes problems when we intend to manually balance Eq.~\eqref{eq:balance1} with generators and DERs from load buses. Specifically, in MATPOWER, $\sum_{k\in \cK} \big(P_k^M - P_k^L(\bm{p}_k,\bm{q}_k)+P_k^0\big)$ is always zero thanks to the slack bus, making updating $\lambda$ with gradient \eqref{eq:gradientLambda} impossible. To address this issue, we apply the following techniques to fix the output of the slack bus.
We record the initial real power of the slack bus as $P^0_{\text{slack}}$. Afterwards, the output of the slack bus changes according to all the other generation and load buses, denoted by $P_{\text{slack}}(\bm{P}^M,\bm{P}^L(\bm{p},\bm{q}))$. Then, instead of enforcing the constraint Eq.~\eqref{eq:balance1}, we use the constraint $P^0_{\text{slack}} = P_{\text{slack}}(\bm{P}^M,\bm{P}^L(\bm{p},\bm{q}))$ to achieve power balance using the generation and load buses while prohibiting the slack bus from participation.

\begin{figure}
	\centering
	\includegraphics[width=.47\textwidth]{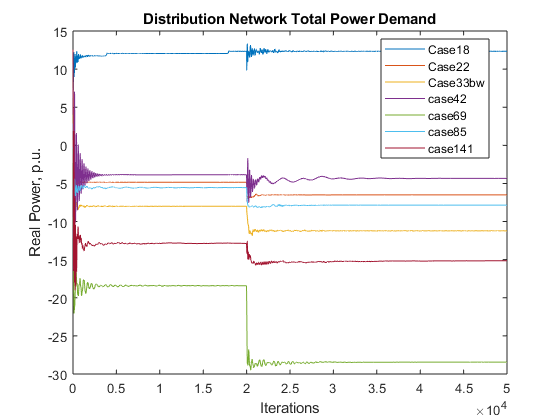}
	\caption{Convergence of total power demands at the substation of all 7 distribution networks.}
	\label{fig:dist_power_converge}
	\vspace{-5mm}
\end{figure}

\begin{figure}
	\centering
	\includegraphics[width=.47\textwidth]{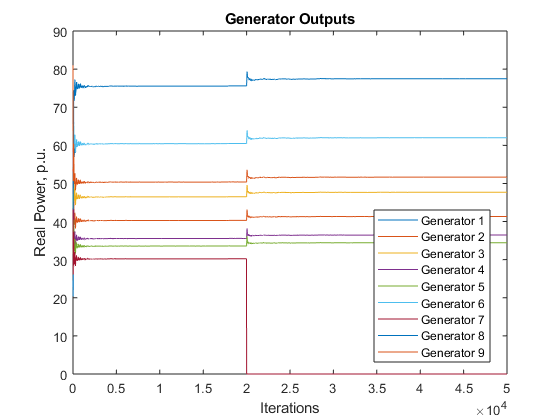}
	\caption{Convergence of power outputs of all 9 generators. Generator 7 goes down after iteration 20,000.}
	\label{fig:gen_converge}
	\vspace{-5mm}
\end{figure}

\begin{figure}
	\centering
	\includegraphics[width=.47\textwidth]{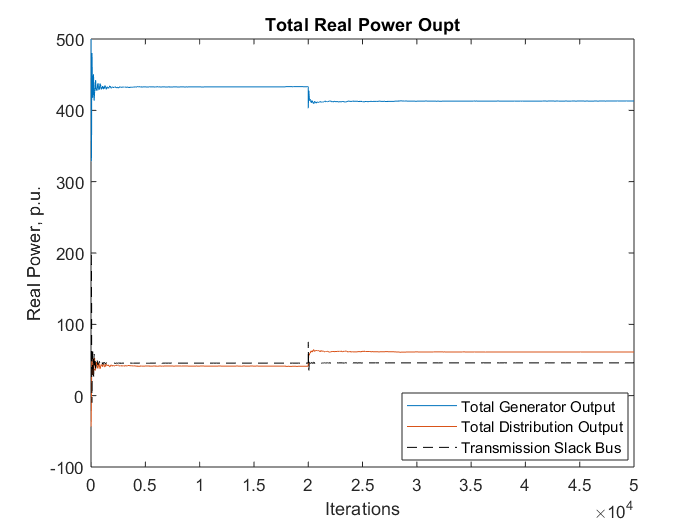}
	\caption{Convergence of total power output from generators and from all distribution networks. The unchanged slack bus output indicates that power demand and supply have been balanced.}
	\label{fig:totaltotal_converge}
	\vspace{-5mm}
\end{figure}

\begin{figure}
	\centering
	\includegraphics[width=.47\textwidth]{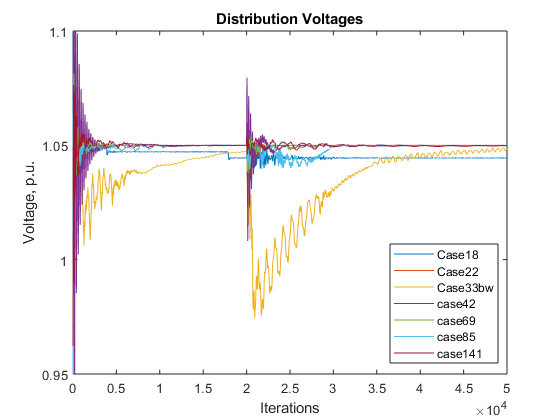}
	\caption{Convergence of voltage magnitudes of sampled nodes from 7 distribution networks. The voltage magnitudes approach the upper bound of 1.05 p.u. to allow for maximal real power generation from DERs.}
	\label{fig:voltage_converge}
	\vspace{-5mm}
\end{figure}

\subsection{Convergence}
\begin{figure}
	\centering
	\includegraphics[width=.47\textwidth]{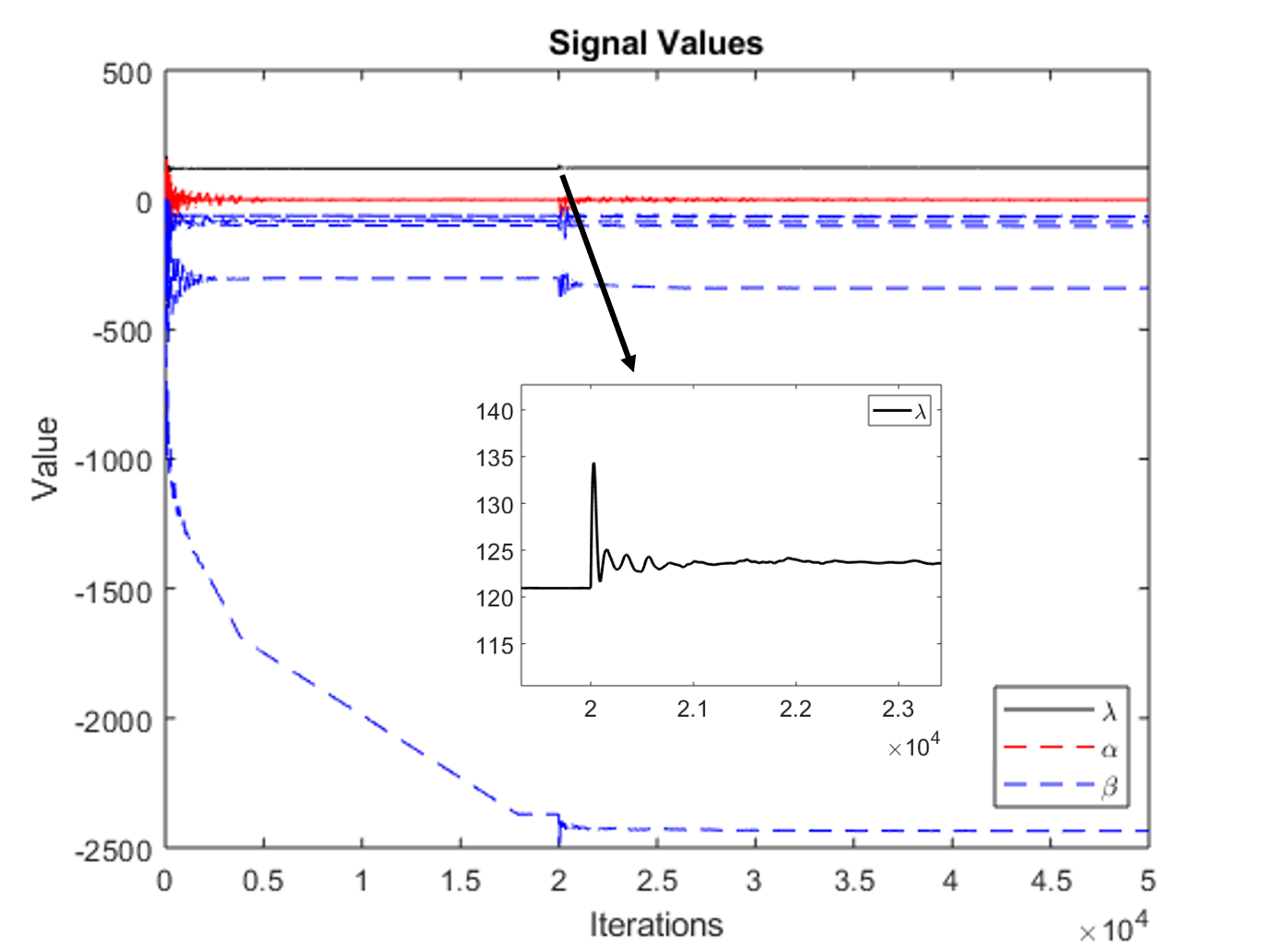}
	\caption{Convergence of $\lambda$, in addition to $\alpha$ and $\beta$ of sampled nodes from all 7 distribution networks.}
	\label{fig:lambda_converge}
	\vspace{-5mm}
\end{figure}

As shown in Fig.~\ref{fig:dist_power_converge} and Fig.~\ref{fig:gen_converge}, the power outputs of the generators and distribution substations gradually stabilize before Generator 7 goes down at iteration 20,000, disturbing the imminent convergence. After Generator 7 goes down, the value of $\lambda$ increases to reflect a higher price for balancing the power demand and supply, as Fig.~\ref{fig:lambda_converge} shows. The larger $\lambda$ then leads to more power outputs from the remaining generators and DERs in the distribution network; see Fig.~\ref{fig:gen_converge} for increasing generator outputs and Fig.~\ref{fig:dist_power_converge} for reducing load bus demand. As a result, we see from Fig.~\ref{fig:totaltotal_converge} that the total outputs from the distribution feeders compensate for the lost generator, with the slack bus maintaining constant output as expected.

Fig.~\ref{fig:voltage_converge} plots the voltage convergence at sampled nodes in all 7 distribution feeders. It is worth noticing that most voltage magnitudes are pushed towards the 1.05~p.u. upper limits because DERs are incentivized to generate more real power to support the transmission system, increasing the voltage levels. Here, negative reactive power is injected to allow for more real power generation without violating the voltage bounds. 

In Fig.~\ref{fig:lambda_converge}, we record and plot the convergence of signal $\lambda$, which is identical for all nodes, along with signals $\alpha$ and $\beta$ of arbitrarily sampled nodes from the 7 distribution feeders without losing generality. Note that the value of $\lambda$ increases after Generator 7 goes down to incentivize more real power input from other generators and DERs. $\alpha$ is smaller than $\lambda$ because $\mathbf{A}_k^{\top}\nabla_{\bm{v}_k} \bm{g}_k(\bm{v}_k)^{\top}\bm{\mu}_k$ in Eq.~\eqref{eq:alpha} is negative to discourage real power generation resulting from overvoltages; $\beta$ is even smaller because we have $\lambda \bm{N}_k=0$ here and $\mathbf{B}_k>\mathbf{A}_k$ element-wise for distribution networks in Eq.~\eqref{eq:beta}, driving more negative reactive power to lower the voltages.

\begin{figure}
	\centering
	\includegraphics[width=.47\textwidth]{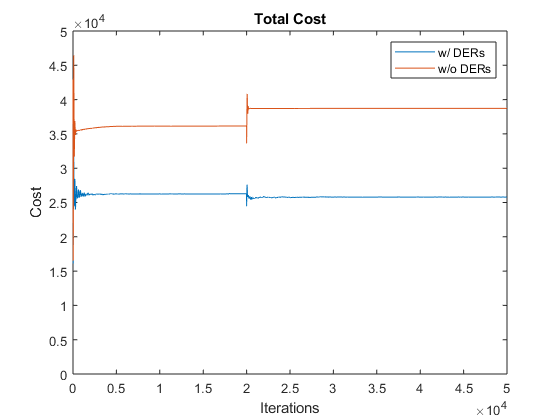}
	\caption{Total cost reduced by replacing some of the high-marginal-cost power from generators with low-marginal-cost power from DERs.}
	\label{fig:totalcost}
	\vspace{-5mm}
\end{figure}

\subsection{Cost Analysis}
We plot the convergence of the values of the total cost function in Fig.~\ref{fig:totalcost}. For comparison, we conduct another set of simulations without involving distribution systems in the economic dispatch, i.e., set $\lambda =0$ in the signals sent to DERs in \eqref{eq:signals}. A higher total cost is recorded without DERs' participation. Such results are expected because even though DERs and generators share similar cost functions, the marginal costs are significantly different: to provide the same amount of power, it would be more economic and more optimal to replace some of the high-marginal-cost generator outputs with low-marginal-cost DERs outputs.

\section{Conclusion}\label{sec:conclude}
This letter formulated a co-optimization problem for economic dispatch in the transmission system and voltage regulation in the distribution networks that are connected to the load buses of the transmission system. Large generators in the transmission system and DERs in the distributed systems are jointly leveraged to achieve the optimal points of the T\&D system. We proposed a primal-dual gradient algorithm, as well as its distributed market-based equivalence, to solve this problem. Numerical results validated that DERs could provide economic grid services to the transmission system while helping to maintain operational constraints in the distribution networks.

\bibliographystyle{IEEEtran}
\bibliography{biblio}

\end{document}